\documentclass[11pt]{article}
\usepackage{amsmath,amsthm}

\def\refeq#1{\if\workingver y(\ref{#1})-[[#1]]\else(\ref{#1})\fi}
\def\refth#1{\if\workingver y\ref{#1}-[[#1]]\else\ref{#1}\fi}
\def\mylabel#1{\if\workingver y\label{#1}{\bf\ \ [[#1]]\ \ }
\else\label{#1}\fi}
\def\mybibitem#1{\if\workingver y\bibitem{#1}{\bf\ \ [[#1]]\ \ }
\else\bibitem{#1}\fi}

\def\institute#1{\gdef\@institute{#1}}

\def\institutename{\par
 \begingroup
 \parskip=\z@
 \parindent=\z@
 \setcounter{@inst}{1}%
 \def\and{\par\stepcounter{@inst}%
 \noindent$^{\the@inst}$\enspace\ignorespaces}%
 \setbox0=\vbox{\def\thanks##1{}\@institute}%
 \ifnum\c@@inst=1\relax
 \else
   \setcounter{footnote}{\c@@inst}%
   \setcounter{@inst}{1}%
   \noindent$^{\the@inst}$\enspace
 \fi
 \ignorespaces
 \@institute\par
 \endgroup}

\newtheorem{thm}{Theorem}
\newtheorem{lem}[thm]{Lemma}
\newtheorem{defn}[thm]{Definition}
\newtheorem{cor}[thm]{Corollary}
\newtheorem{prop}[thm]{Proposition}

\newtheorem{rem}[thm]{Remark}
\newtheorem{conj}[thm]{Conjecture}
\newtheorem{example}[thm]{Example}

\def\begeq#1{\begin{equation}\mylabel{#1}}
\def\endeq{\end{equation}}

\def\begalg{\begin{alg}}
\def\endalg{\end{alg}}

\let\workingver=n

\newcommand{\beq}{\begin{equation}}

\font\tenrm=cmr10

\def\ga{\gamma}
\def\de{\delta}
\def\al{\alpha}
\def\be{\beta}
\begin{document}

\title{Netted Binomial Matrices}
\author{
\vspace{1cm}
Pantelimon St\u anic\u a\thanks{On leave from the Institute of
Mathematics of Romanian Academy, Bucharest, Romania}\\
\tenrm Auburn University Montgomery\\
\tenrm Department of Mathematics,\\
\tenrm Montgomery, AL 36117, USA\\
\tenrm e-mail: stanpan@strudel.aum.edu
}
\date{\today}
\maketitle

%\pagestyle{myheadings}

%\begin{abstract}
%We prove lots of stuff.
%\end{abstract}

\baselineskip=1.98\baselineskip

\newpage

{\em Running Head:}\\
{\bf Netted Binomial Matrices}

{\em Address for Proofs:}

\noindent
Pantelimon Stanica\\
Department of Mathematics\\
Auburn University Montgomery\\
Montgomery, AL 36117\\
USA

\newpage
{\bf Abstract:}

{\em
We prove that powers of $3$-netted matrices (the entries
satisfy a third-order recurrence
$\de a_{i,j}=\al a_{i-1,j}+\be a_{i-1,j-1} +\ga a_{i,j-1}$) preserve
this property of nettedness, that is the entries of
the $e$-th power satisfy
$\de_{e} a_{i,j}^{(e)}=\al_e a_{i-1,j}^{(e)}+
\be_e a_{i-1,j-1}^{(e)}+\ga_e a_{i,j-1}^{(e)},$
where the coefficients are all instances of the same sequence,
$x_{e+1}=(\be+\de)x_e-(\be\de+\al\ga) x_{e-1}$.
 Also, we find a matrix $T_n(m)$
 and a vector $v$, such that
$T_n(m)^e\cdot v$ gives $n$ consecutive entries of the general Fibonacci
(Pell) sequence with parameter $m$.
It generalizes the known property $\displaystyle
\left(\begin{array}{cc}
 0&1\\
 1&1
  \end{array}\right)^e\cdot (1,0)^t=(F_{e-1},F_{e})^t.
$
We close by giving a conjecture about the spectral
properties of the binomial matrix we found.}

\newpage

\section{Introduction}

In \cite{PS}, Peele and St\u anic\u a studied $n\times n$
matrices with the $(i,j)$
 entry the binomial coefficient $\binom{i-1}{j-1}$
 (matrix $L_n$), respectively
$\binom{i-1}{n-j}$ (matrix $R_n$)
 and derived many interesting results
on powers of these matrices.
$L_n$ was easily subdued, but
curiously enough, closed forms for entries of powers of $R_n$, say $R_n^e$,
were not found. However,  recurrences among various entries of
 $R_n^e$ were
proved and precise results on congruences modulo any prime $p$
 were found.
They proved that the entries $a_{i,j}^{(e)}$ of the $e$-th
power of $R_n$, satisfy
\[
F_{e-1} a_{i,j}^{(e)}=F_{e} a_{i-1,j}^{(e)}+
F_{e+1} a_{i-1,j-1}^{(e)}-F_{e} a_{i,j-1}^{(e)},
\]
where $F_e$ is the Fibonacci sequence, $F_{e+1}=F_e+F_{e-1},
F_0=0,F_1=1$,
which was used to prove beautiful results on powers of these matrices
modulo a prime $p$.
As we shall see in our second section,
this is not a singular phenomenon.
In this paper we generalize the results of \cite{PS}
for a class of matrices, containing $R_n$,
where the entries satisfy a certain
 recurrence (we call these {\em netted} matrices).
Moreover, we find a matrix with the property that any power
multiplied by a fixed vector
gives a tuple of consecutive terms of the Pell or
Fibonacci numbers sequence.
We also find the generating function for the entries of powers
of these matrices.
As applications we find four interesting identities for
Fibonacci (or Pell) numbers.
In the fifth section we provide a few results on the order
of these matrices
modulo a prime and in the last section we propose a
conjecture on the spectral properties of $T_n(m)$.

\section{Sequences Satisfying a Third Order Recurrence }

 Define a
tableau with elements $a_{i,j},\, i\geq 0,\, j\geq 0$, which satisfy
(for $i\geq 1,j\geq 1$)
\begin{equation}
\label{(1)}
\de a_{i,j}=\al a_{i-1,j}+\be a_{i-1,j-1} +\ga a_{i,j-1},
\end{equation}
with the boundary conditions
\begin{eqnarray}
&& \be a_{i,0}+\ga a_{i+1,0}=0, \forall\ 1\leq i\leq n-1 \label{bound1}\\
&& \de a_{i+1,n+1}-\al a_{i,n+1}=0, \forall\ 1\leq i\leq n-1 \label{bound2}.
\end{eqnarray}
We remark that if the $0$-th and $(n+1)$-th column are made up of
zeros, then the conditions \refeq{bound1} and \refeq{bound2} are
fulfilled.

 In our main result of this section we prove that  \refeq{(1)} is
 preserved for higher powers of the $n\times n$ matrix
 $(a_{i,j})_{i=1\ldots n,j=1\ldots n}$.
 Precisely, we prove
\begin{thm}
\label{general_theorem}
The entries of the $e$-th power of the matrix
$R=(a_{i,j})_{i=1...n,j=1...n}$ satisfy the recurrence
\[
\de_{e} a_{i,j}^{(e)}=\al_e a_{i-1,j}^{(e)}+
\be_e a_{i-1,j-1}^{(e)}+\ga_e a_{i,j-1}^{(e)},
\]
where the sequences $\al_e,\be_e,\ga_e,\de_e$ are all instances of the
 sequence $x_e$ satisfying
\[
x_{e+1}=(\be+\de)x_e-(\be\de+\al\ga) x_{e-1},
\]
with initial conditions
$(\de_1 =\de;  \de_2=\de^2-\al\ga)$;
$(\al_1  =\al;   \al_2=\al(\de+\be))$;
$(\be_1 =\be;  \be_2=\be^2-\al\ga)$ and
$(\ga_1= \ga;  \ga_2=\ga(\be+\de))$.
\end{thm}
\begin{proof}
We prove by induction on $e$ that there exists a relation among
 the entries of any $2\times 2$ cells, namely
\[
\delta_e a_{i,j}^{(e)}=\alpha_e a_{i-1,j}^{(e)}+\beta_e
a_{i-1,j-1}^{(e)}+\gamma_e a_{i,j-1}^{(e)}.
\]
The above relation is certainly true for $n=1$.
We evaluate, for $i\geq 2$,
\begin{eqnarray*}
\al\de_{e-1} a_{i-1,j}^{(e)}
&= &
 \sum_{s=1}^n \al\de_{e-1} a_{i-1,s}  a_{s,j}^{(e-1)}\\
&=&
\sum_{s=1}^n \al a_{i-1,s}  \left(\al_{e-1} a_{s-1,j}^{(e-1)}+
\be_{e-1}a_{s-1,j-1}^{(e-1)}+
\ga_{e-1} a_{s,j-1}^{(e-1)}\right)\\
&= &
\sum_{s=1}^n\left (\de a_{i,s} -\be a_{i-1,s-1} -\ga a_{i,s-1} \right)
\left(\al_{e-1}
a_{s-1,j}^{(e-1)} +\be_{e-1} a_{s-1,j-1}^{(e-1)}\right)  \\
&& +
 \sum_{s=1}^n \al\ga_{e-1} a_{i-1,s} a_{s,j-1}^{(e-1)}=
 \sum_{s=1}^n \de a_{i,s}
\left (\al_{e-1} a_{s-1,j}^{(e-1)}+\be_{e-1}  a_{s-1,j-1}^{(e-1)}\right)
\\
%\end{eqnarray*}
%\begin{eqnarray*}
&& -
\ga\al_{e-1}  a_{i,j}^{(e)} - \be\be_{e-1}  a_{i-1,j-1}^{(e)} -
\ga\be_{e-1} a_{i,j-1}^{(e)}-\be\al_{e-1}  a_{i-1,j}^{(e)}\\
&&
+\al\ga_{e-1}  a_{i-1,j-1}^{(e)}  -
\ga\al_{e-1}\left( a_{i,0}a_{0,j}^{(e-1)}-a_{i,n}a_{n,j}^{(e-1)} \right)\\
&& -
\be\be_{e-1}\left( a_{i-1,0}a_{0,j-1}^{(e-1)}-a_{i-1,n}a_{n,j-1}^{(e-1)}
\right)\\
&& -
\ga\be_{e-1}\left( a_{i,0}a_{0,j-1}^{(e-1)}-a_{i,n}a_{n,j-1}^{(e-1)}
\right)\\
&& -
\be\al_{e-1}\left( a_{i-1,0}a_{0,j}^{(e-1)}-
a_{i-1,n}a_{n,j}^{(e-1)} \right).
\end{eqnarray*}
Using the boundary conditions \refeq{bound1}
 and \refeq{bound2},  we obtain, for $i\geq 2$,
\begin{eqnarray*}
\al\de_{e-1}a_{i-1,j}^{(e)}
&=&
(\al \ga_{e-1}-\be \be_{e-1})  a_{i-1,j-1}^{(e)}-
\ga\al_{e-1}  a_{i,j}^{(e)}- \ga\be_{e-1} a_{i,j-1}^{(e)}\\
&& -
 \be\al_{e-1}  a_{i-1,j}^{(e)}+
\sum_{s=1}^n \de a_{i,s}  \left (\de_{e-1}  a_{s,j}^{(e-1)}-
\ga_{e-1}  a_{s,j-1}^{(e-1)}\right ) \\
&& +
\left( \al_{e-1} a_{n,j}^{(e-1)}+\be_{e-1} a_{n,j-1}^{(e-1)} \right)
\left(\be a_{i-1,n}+\ga a_{i,n}  \right)\\
&& -
\left( \al_{e-1} a_{0,j}^{(e-1)}+\be_{e-1} a_{0,j-1}^{(e-1)} \right)
\left(\be a_{i-1,0}+\ga a_{i,0}  \right)
\\
%\end{eqnarray*}
%\begin{eqnarray*}
&=&
(\al \ga_{e-1}-\be \be_{e-1})  a_{i-1,j-1}^{(e)}-
\ga\al_{e-1}  a_{i,j}^{(e)}-\ga\be_{e-1} a_{i,j-1}^{(e)}\\
&& - \be\al_{e-1}  a_{i-1,j}^{(e)}+
\de \de_{e-1}  a_{i,j}^{(e)}-\de\ga_{e-1}  a_{i,j-1}^{(e)}\\
&& +
\left( \al_{e-1} a_{n,j}^{(e-1)}+\be_{e-1} a_{n,j-1}^{(e-1)} \right)
\left(\de a_{i,n+1}-\al a_{i-1,n+1}  \right)
%\\
%\\
\end{eqnarray*}
\begin{eqnarray*}
&=&
(\al \ga_{e-1}-\be \be_{e-1})  a_{i-1,j-1}^{(e)}+
(\de\de_{e-1}-\ga\al_{e-1})  a_{i,j}^{(e)}\\
&& -
(\ga\be_{e-1}+\de\ga_{e-1}) a_{i,j-1}^{(e)}
- \be\al_{e-1}  a_{i-1,j}^{(e)}.
\end{eqnarray*}
Thus,
\[
\begin{split}
(\de\de_{e-1}-\ga\al_{e-1}) a_{i,j}^{(e)} =(\al\de_{e-1}+  \be\al_{e-1})
a_{i-1,j}^{(e)} & +
(\be\be_{e-1}-\al \ga_{e-1})  a_{i-1,j-1}^{(e)} \\
&  + (\ga\be_{e-1}+\de\ga_{e-1}) a_{i,j-1}^{(e)}.
\end{split}
\]
Therefore, we obtain the system of sequences
\begin{eqnarray}
\de_e  & =& \de\de_{e-1}-\ga\al_{e-1}\label{eq14}\\
\al_e   & =& \al\de_{e-1}+\be \al_{e-1}\label{eq15}\\
\be_e & = & \be\be_{e-1}-\al\ga_{e-1}\label{eq16}\\
\ga_e & = & \ga\be_{e-1}+\de\ga_{e-1}.\label{eq17}
\end{eqnarray}
From \refeq{eq14} we get $\al_{e-1}=(\de/\ga) \de_{e-1}-(1/\ga) \de_e$,
 which
replaced in \refeq{eq15} gives the recurrence
\[
\de_{e+1}=(\be+\ga) \de_e-(\be\de+\ga\al)\de_{e-1}.
\]
Similarly,
\begin{eqnarray*}
\al_{e+1}&=&(\be+\ga) \al_e-(\be\de+\ga\al)\al_{e-1}\\
\be_{e+1}&=&(\be+\ga) \be_e-(\be\de+\ga\al)\be_{e-1}\\
\ga_{e+1}&=&(\be+\ga) \ga_e-(\be\de+\ga\al)\ga_{e-1}.
\end{eqnarray*}
The initial conditions are
$(
\de_1 =\de;  \de_2=\de^2-\al\ga),
(\al_1 =\al;   \al_2=\al(\de+\be)),
(\be_1 =\be;  \be_2=\be^2-\al\ga),
(\ga_1= \ga;  \ga_2=\ga(\be+\de))
$.
\end{proof}

\begin{example}
As examples of tableaux satisfying our conditions, we have
$a_{i,j}^1=\binom{i-1}{j-1}$ $(\de=1,\al=1,\be=1,\ga=0)$,
$a_{i,j}^2=\binom{i-1}{n-j}$ $(\de=0,\al=1,\be=1,\ga=-1)$,
$a_{i,j}^3=\binom{n-i}{n-j}$ $(\de=1,\al=0,\be=-1,\ga=1)$.
Other examples are given by the alternating matrices $(-1)^{i+j} a_{i,j}^k$
(or  $(-1)^{i-1} a_{i,j}^k$ or $(-1)^{j-1} a_{i,j}^k$, etc.), $k=1,2,3$.
In the next section we present more examples.
\end{example}

\section{Fibonacci and Pell Matrices}

In this section we uncover a very interesting side
of the previous section's results.
A matrix of the form
$M=
\left(\begin{array}{cc}
0 & 1 \\
1 & m
\end{array}\right)
$
is called a {\em Fibonacci matrix}. It is known that if
the sequence
 $U_{e+1}=mU_{e}+U_{e-1},\ U_0=0,\, U_1=1$, then
$M^e=
\left(\begin{array}{cc}
U_{e-1} & U_e \\
U_e & U_{e+1}
\end{array}\right)
$
and $M^e\cdot
\left(\begin{array}{c}
1  \\
0
\end{array}\right)=
\left(\begin{array}{c}
U_{e-1}  \\
U_e
\end{array}\right).
$
If $m$ is an indeterminate, then $U_e$ is called the
{\em Fibonacci polynomial}.
If $m=2$, $U_e (=P_e)$ is the {\em Pell sequence}.
The question that arises is whether there are higher
dimensional square
 matrices $T_n(m)$ such that $T_n(m)^e\cdot v$ is a vector of $n$
 consecutive terms of
the sequence $U_n$, for some vector $v$ and any power $e$.
We are able to answer positively the posed question
for such a sequence. Let $I_n$ be
the identity matrix of dimension $n$ and
 $M^t$ be the {\em transpose} of a given matrix $M$.

First, we consider the Pell sequence:
$P_{e+1}=2 P_e+P_{e-1},\, P_0=0,\, P_1=1$.
We prove that for each dimension there is a unique
$n\times n$ matrix $T_n$
 with positive entries, constructed by bordering $T_{n-1}$
 and such that
 the entries satisfy $\delta a_{i,j}=\alpha a_{i-1,j}+
 \beta a_{i-1,j-1}+
 \gamma a_{i,j-1}$, with $\delta =0, \alpha=1,\beta=2,\gamma=-1$.
 We also
 prove, using the previous section's result, that the entries of
  any power of this matrix satisfy a similar relationship,
   where the corresponding coefficients are all instances of
   the {\em Pell sequence}.
Let $a_{i,j}= a_{i,j}^{(1)}=2^{i+j-n-1} \binom{i-1}{n-j},\ i,j\geq 0$.

\begin{thm}
\label{pell1}
Let $v=((-1)^n P_{n-1},(-1)^{n-1} P_{n-2},\ldots,-P_0)^t$ and
$T_n=(a_{i,j})_{1\leq i,j\leq n}$.
We have $T_n^{e+1} \cdot v=(P_{(n-1)e},P_{(n-1)e+1},
\ldots,P_{(n-1)(e+1)})^t$ and
\[
P_{e-1} a_{i,j}^{(e)}+P_e a_{i,j-1}^{(e)}=P_e a_{i-1,j}^{(e)}+P_{e+1}
a_{i-1,j-1}^{(e)},
\]
where $a_{i,j}^{(e)}$ are the entries of $T_n^e$ and $P_e$ is the
Pell sequence. Moreover, $T_n$ is unique with the property
$a_{1,j}=0, j<n, a_{i,n}=2^{i-1}$ and $a_{i,j}=2 a_{i-1,j}+a_{i-1,j+1}$.
\end{thm}

\begin{proof}
By induction on $e$ we prove that
$
\sum_{j=1}^n (-1)^{n+1-j} a_{i,j}^{(e+1)} P_{n-j}=P_{(n-1)e+i-1},
$
which will imply the first assertion.
Assume $e=0$. We need to show
$\sum_{j=1}^n (-1)^{n+1-j} a_{i,j} P_{n-j}=P_{i-1}$, which will be proved
by showing that the left hand side expression satisfies the Pell recurrence
with the initial conditions of $P_{i-1}$. Denote by $X_{i-1}$ the left hand
side expression. First,
$X_0=\sum_{j=1}^n (-1)^{n+1-j} a_{1,j} P_{n-j}=
(-1)^{n+1-n} a_{1,n} P_{n-n}=0$.
Now,
$X_1=\sum_{j=1}^n (-1)^{n+1-j} a_{2,j} P_{n-j}=
(-1)^{1} a_{2,n} P_{0}+(-1)^{2} a_{2,n} P_{1}=1$.
Assume $1\leq i\leq n-2$. Then
\begin{eqnarray*}
X_{i+1}
&=& \sum_{j=1}^n (-1)^{n+1-j} a_{i+2,j} P_{n-j}\\
&=&
\sum_{j=1}^n (-1)^{n+1-j} (2 a_{i+1,j}+a_{i+1,j+1} ) P_{n-j}\\
&=& 2 X_i +\sum_{j=1}^n (-1)^{n+1-j} a_{i+1,j+1} P_{n-j}\\
%\end{eqnarray*}
%\begin{eqnarray*}
& =& 2 X_i+
\sum_{j=1}^n (-1)^{n+1-j} a_{i+1,j+1} (2 P_{n-j-1}+P_{n-j-2})\\
&= & 2X_i+
2\sum_{j=1}^{n-1} (-1)^{n+1-j} a_{i+1,j+1} P_{n-(j+1)}\\
& &+
\sum_{j=1}^{n-1} (-1)^{n+1-j} a_{i+1,j+1} P_{n-(j+1)-1}\\
&\stackrel{j+1=s}{=} & 2X_i+
2 \sum_{s=2}^{n} (-1)^{n+2-s} a_{i+1,s} P_{n-s}+
\sum_{s=2}^{n} (-1)^{n+2-s} a_{i+1,s} P_{n-s-1}\\
&=& 2 (-1)^{n} a_{i+1,1} P_{n-1}-Y_{i} +(-1)^n a_{i+1,1} P_{n-2}\\
&=&
 (-1)^{n} a_{i+1,1} P_{n}-Y_i
 \stackrel{a_{i+1,1}=0\ \text{if}\ i\leq n-2}{=} -Y_i,
\end{eqnarray*}
where
\begin{eqnarray*}
Y_{i}&\stackrel{def}{=}&\sum_{s=1}^{n} (-1)^{n+1-s} a_{i+1,s} P_{n-s-1}\\
&=& \sum_{s=1}^{n} (-1)^{n+1-s} (2a_{i,s}+a_{i,s+1}) P_{n-s-1}\\
&=& 2 \sum_{s=1}^{n} (-1)^{n+1-s} a_{i,s} P_{n-s-1}+
 \sum_{s=1}^{n} (-1)^{n+1-s} a_{i,s+1} P_{n-s-1}\\
 &=& 2 Y_{i-1}+ \sum_{u=2}^{n} (-1)^{n+2-u} a_{i,u} P_{n-u}\\
 &=& 2 Y_{i-1}-X_{i-2} +(-1)^n a_{i,1} P_{n-1}
 \stackrel{a_{i,1}=0\ \text{if}\ i\leq n-2}{=}
 2 Y_{i-1}-X_{i-2}
\end{eqnarray*}
Using $Y_i=-X_{i+1}$ in the previous
recurrence we get
\[
X_{i+1}=2 X_{i} +X_{i-1},
\]
relation satisfied by the Pell sequence. Since $X_{i-1}$
has the same initial
conditions as $P_{i-1}$ we have
 $\sum_{j=1}^n (-1)^{n+1-j} a_{i,j} P_{n-j}=P_{i-1}$.

 The first step of induction is proven. Now,
 \[\begin{split}
&
\sum_{j=1}^n (-1)^{n+1-j} a_{i,j}^{(e+1)} P_{n-j}=
\sum_{j=1}^n (-1)^{n+1-j} \sum_{k=1}^na_{i,k} a_{k,j}^{(e+1)} P_{n-j}\\
&
=\sum_{k=1}^n a_{i,k} \sum_{j=1}^n (-1)^{n+1-j} a_{k,j}^{(e)} P_{n-j}\\
&
= \sum_{j=1}^n a_{i,k} P_{(n-1)(e-1)+k-1}.
\end{split}
 \]
 We shall prove that the matrix
 $T$ acts as an index-translation on the Pell sequence, namely
 \[
\sum_{k=1}^n a_{i,k} P_{t+k}=P_{t+n+i-1},\ t\geq -1.
 \]
 If this is so, then
by taking $t=(n-1)(e-1)-1$, the step of induction will be done.
Let $W_i=\sum_{k=1}^n a_{i,k} P_{t+k}$ ($t$ is assumed fixed).
First, $W_1=\sum_{k=1}^n a_{1,k} P_{t+k}=a_{1,n}P_{t+n}=P_{t+n}$.
Then, $W_2=\sum_{k=1}^n a_{2,k} P_{t+k}=a_{2,n-1}P_{t+n-1}+a_{2,n}P_{t+n}=
P_{t+n-1}+2P_{t+n}=P_{t+n+1}$.
Now, for $1\leq i\leq n-1$,
\begin{eqnarray*}
W_{i+1}
&=&\sum_{k=1}^n a_{i+1,k} P_{t+k} =
\sum_{k=1}^n (2 a_{i,k}+a_{i,k+1}) P_{t+k}\\
&=& 2W_i+\sum_{k=1}^n a_{i,k+1} P_{t+k}=2 W_i+\sum_{k=1}^{n-1}
a_{i,k+1}(P_{t+k+2}-2 P_{t+k+1})\\
&\stackrel{u=k+1}{=}&
2W_i+\sum_{u=2}^n a_{i,u} P_{t+u+1}-2 \sum_{u=2}^n a_{i,u} P_{t+u}\\
&= &
V_i-a_{i,1} P_{t+2}+2 a_{i,1} P_{t+1}
\stackrel{a_{i,1}=0\ \text{if}\ i\leq n-1}{=}
V_i,
\end{eqnarray*}
where
\begin{eqnarray*}
V_i&=&\sum_{u=1}^n a_{i,u} P_{t+u+1}=
\sum_{u=1}^n (2 a_{i-1,u}+a_{i-1,u+1}) P_{t+u+1}\\
&\stackrel{a_{i+1,n+1}=0}{=}&
2 V_{i-1}+\sum_{u=1}^{n-1} a_{i-1,u+1} P_{t+u+1}
\stackrel{u+1=s}{=} 2 V_{i-1} +\sum_{u=2}^n a_{i-1,s} P_{t+s}\\
&=& 2 V_{i-1}+X_{i-1}-a_{i-1,1} P_{t+1}
\stackrel{a_{i-1,1}=0\ \text{if}\ i\leq n-1}{=} 2V_{i-1}+W_{i-1}.
\end{eqnarray*}
Using $V_i=W_{i+1}$ in the previous recurrence, we get
$W_{i+1}=2 W_i +W_{i-1}$. Therefore, $W_i=P_{t+n+i-1}$, since
$W_1=P_{t+n}, W_2=P_{t+n+1}$.

Using Theorem \refth{general_theorem}, with
$\delta =0, \alpha=1,\beta=2,\gamma=-1$, we get the recurrence between
the entries of the higher power of $T_n$, namely
\(
P_{e-1} a_{i,j}^{(e)}+P_e a_{i,j-1}^{(e)}=P_e a_{i-1,j}^{(e)}+P_{e+1}
a_{i-1,j-1}^{(e)}.
\)

The fact that $T_n$ is the unique matrix with the given properties
follows easily observing that such a matrix could be
defined inductively as follows:
let $T_1=1$. Assume
$T_{n-1}=\left( a_{i,j}\right)_{i,j=1,2,\ldots,n-1}$ and
construct $T_n$ by bordering
$T_{n-1}$ with the first column and the last row (left and bottom).
The first column is $(0,0,\ldots,0,1)^t$ and the last row is given by:
$a_{n,n}=2^{n-1}$ and $a_{n,j}=2 a_{n-1,j}+a_{n-1,j+1}$.
\end{proof}

Let $a_{i,j}=a_{i,j}^{(1)}=m^{i+j-n-1} \binom{i-1}{n-j}$.
Similarly, we can show (we omit the proof)
\begin{thm}
\label{general_fibonacci}
Let $w=((-1)^n U_{n-1},(-1)^{n-1} U_{n-2},\ldots,-U_0)^t$ and
$T_n(m)=(a_{i,j})_{i,j}$.
Then $T_n(m)^{e+1} \cdot w=(U_{(n-1)e},U_{(n-1)e+1},
\ldots,U_{(n-1)(e+1)})^t$ and
\begin{equation}
\label{gen_fib}
U_{e-1} a_{i,j}^{(e)}+U_e a_{i,j-1}^{(e)}=U_e a_{i-1,j}^{(e)}+U_{e+1}
a_{i-1,j-1}^{(e)},
\end{equation}
where $a_{i,j}^{(e)}$ are the entries of $T_n(m)^e$ and $U_e$ is the
sequence satisfying $U_{e+1}=m U_e+U_{e-1},\ U_0=0,\, U_1=1$.
 Moreover, $T_n(m)$ is unique with the property
$a_{1,j}=0, j<n, a_{i,n}=m^{i-1}$ and $a_{i,j}=m a_{i-1,j}+a_{i-1,j+1}$.
\end{thm}

\begin{defn}
We call such a matrix $T_n(m)$ a {\em generalized Fibonacci
matrix} of dimension $n$
and parameter $m$. If $m=2$, $T_n(2)=T_n$ is the {\em Pell matrix}.
\end{defn}

\begin{example}
We give here the first few powers of $T_3(m)$,
\[
\begin{split}
&
T_3(m)=
\left(
\begin{matrix} 
0 & 0 & 1 \cr 0 & 1 & m \cr 1 & 2\,m & m^2 \cr  
\end{matrix}
\right),\ 
T_3(m)^2=
\left(
  \begin{matrix}
 1 & 2\,m & m^2 \cr m & 1 + 2\,m^2 & m + m^3 \cr m^2 & 2\,
    \left( m + m^3 \right)  & {\left( 1 + m^2 \right) }^2 \cr  
\end{matrix}
\right),\\
&
T_3(m)^3=
\left(
  \begin{matrix}
 m^2 & 2\,\left( m + m^3 \right)  & {\left( 1 + m^2 \right) }^2 \cr m + 
    m^3 & 1 + 4\,m^2 + 2\,m^4 & m\,\left( 2 + 3\,m^2 + m^4 \right)  \cr {\left(
       1 + m^2 \right) }^2 & 2\,m\,\left( 2 + 3\,m^2 + m^4 \right)  & m^2\,
    {\left( 2 + m^2 \right) }^2 \cr  
\end{matrix}
\right) 
\end{split}
\]
\end{example}

By taking some particular cases of our previous results we get
some very interesting binomial sums.
For instance,
\begin{cor}
We have
\item[1.]
\(\displaystyle
\sum_{j=1}^n  (-1)^{n+1-j}  m^{i+j-n-1}\binom{i-1}{n-j}
U_{n-j}=U_{i-1}.
\)
\item[2.]
\(\displaystyle
\sum_{j=1}^n \sum_{k=1}^n (-1)^{n+1-j}m^{i+j+2k-2n-2}
\binom{i-1}{n-k}\binom{k-1}{n-j}U_{n-j}=U_{n+i-2}.
\)
\item[3.]
$\displaystyle
\sum_{j=1}^n U_{l-1}^{n-j}U_l^{j-1} U_{(n-1)p+j-1}\binom{n-1}{j-1}=
U_{(n-1)(l+p)}$, for any $l,p.$
\item[4.]
$\displaystyle
\sum_{j=1}^n U_{(n-1)p+j-1} U_{l-1}^{n-j-1} U_l^{j-2}
\left[ U_l^2 \binom{n-1}{j-1}+
(-1)^l \binom{n-2}{j-2} \right]=
U_{(n-1)(l+p)+1}$, for any $l,p$.
\end{cor}
\begin{proof}
Using Theorem \refth{general_fibonacci}, with $e=1,2$,
we obtain the first two identities.
Now, with the help of Theorem \refth{general_fibonacci} and the
 trivial identity
 $T_n(m)^{l+p}=T_n(m)^l T_n(m)^p$, we get
\begin{eqnarray*}
\left(T_n(m)^l T_n(m)^p\right)\cdot v
&=&
T_n(m)^l\cdot \left(U_{(n-1)p},U_{(n-1)p+1},\ldots,
U_{(n-1)(p+1)}\right)^t\\
&=&
\left(U_{(n-1)(l+p)},U_{(n-1)(l+p)+1},\ldots,U_{(n-1)(l+p+1)}\right).
\end{eqnarray*}
Since $a_{1,j}^{(l)} = U_{l-1}^{n-j} U_l^{j-1}\binom{n-1}{j-1}$,
we obtain the third identity.
Using $a_{2,j}^{(l)}=U_{l-1}^{n-j-1} U_l^{j}\binom{n-2}{j-1}+
U_{l-1}^{n-j} U_e^{j-2}U_{l+1}\binom{n-2}{j-2}$ and
Cassini's identity (see \cite{GKP}, p. 292) (usually given
for the Fibonacci numbers, but certainly true for
the sequence $U_l$, as well, as the reader can check easily,
by induction)
\(
U_{l-1} U_{l+1}-U_l^2=(-1)^l,
\)
we get the fourth identity.
\end{proof}

 \begin{rem}
In general, we have
\[
\displaystyle
\sum_{j=1}^n U_{(n-1)p+j-1} a_{i,j}^{(l)}=U_{(n-1)(l+p)+i-1}
\]
for any $i,l,p$.
 \end{rem}

%Using the following formulas for the Pell sequence,
%\(
%P_e=\sum_{k=1}^e 2^{e-2k+1} \binom{e-k}{k-1},
%\)
%or, in general
%$U_e=\sum_{k=1}^e m^{e-2k+1} \binom{e-k}{k-1}$
%we can write  the previous corollary
 %in terms of binomial coefficients (the last one becomes very
 %complicated).

In general, finding closed forms for the entries of powers of $T_n(m)$
seems to be a very difficult matter.
We can derive (after some work) simple formulas for the
entries of the second row and column of $T_n(m)^e$.
\begin{prop}
We have
$a_{2,j}^{(e)}=U_{e-1}^{n-j-1} U_e^{j}\binom{n-2}{j-1}+
U_{e-1}^{n-j} U_e^{j-2}U_{e+1}\binom{n-2}{j-2}$, and
$a_{i,2}^{(e)}=(n-i)U_{e-1}^{n-i} U_e^{i-1}+
(i-1) U_{e-1}^{n-i+1} U_e^{i-2} U_{e+1}.$
\end{prop}
\begin{rem}
Since $a_{j,n}^{(e-1)}=a_{j,1}^{(e)},
a_{n,j}^{(e-1)}=a_{1,j}^{(e)}$, we get closed forms for
the last row and column of $T_n(m)^e$, as well.
\end{rem}

\newpage

\section{Some Generating Functions and an Inverse}

Although we cannot find simple closed forms for {\em all} entries of
$T_n(m)^e$, we prove
\begin{thm}
The generating function for $a_{i,j}^{(e)}$ is
\[
\displaystyle
B_n^{(e)}(x,y)=
\frac{(U_{e-1}+U_e y)^n}{ U_{e-1}+U_e y-x( U_e+U_{e+1} y)}
\]
\end{thm}

\begin{proof}
Multiplying the recurrence \refeq{gen_fib} by
$x^{i-1} y^{j-1}$ and summing for $i,j\geq 2$, we get
\begin{eqnarray*}
&&U_{e-1}\sum_{i,j\geq 2} a_{i,j}^{(e)} x^{i-1} y^{j-1} +
U_e y \sum_{i,j\geq 2} a_{i,j-1}^{(e)} x^{i-1} y^{j-2}=\\
&&
U_{e}x\sum_{i,j\geq 2} a_{i-1,j}^{(e)} x^{i-2} y^{j-1} +
U_{e+1}xy\sum_{i,j\geq 2} a_{i-1,j-1}^{(e)} x^{i-2} y^{j-2}
\end{eqnarray*}
Thus,
\begin{eqnarray*}
&&
U_{e-1}
\left(
B_n^{(e)}(x,y)-\sum_{i\geq 1} a_{i,1}^{(e)} x^{i-1}-\sum_{j\geq 1}
a_{1,j}^{(e)} y^{j-1}+a_{1,1}^{(e)}
 \right)
 \\
 &&\hspace{2.6cm} +
 U_e y\left( B_n^{(e)}(x,y)-\sum_{j\geq 1} a_{1,j}^{(e)} y^{j-1} \right)=\\
&&
U_e x\left(B_n^{(e)}(x,y)-\sum_{i\geq 1} a_{i,1}^{(e)} x^{i-1}\right)+
U_{e+1} x y B_n^{(e)}(x,y).
\end{eqnarray*}
Solving for $B_n^{(e)}(x,y)$, we get
\begin{equation}
\label{first_eq_gen}
\begin{split}
&
B_n^{(e)}(x,y) \left(U_{e-1}+U_e y-x( U_e+U_{e+1} y)  \right)=\\
&
(U_{e-1}-U_e x) \sum_{i\geq 1} a_{i,1}^{(e)} x^{i-1}+
(U_{e-1}+U_e y) \sum_{j\geq 1} a_{1,j}^{(e)} y^{j-1}-U_{e-1}a_{1,1}^{(e)}.
\end{split}
\end{equation}
We need to find $a_{i,1}^{(e)}$ and $a_{1,j}^{(e)}$.
We prove
\begin{equation}
\label{first_row}
\begin{split}
a_{1,j}^{(e)} &=\binom{n-1}{j-1} U_{e-1}^{n-j} U_e^{j-1}\\
a_{i,1}^{(e)} &=U_{e-1}^{n-i} U_e^{i-1}.
\end{split}
\end{equation}
There is no difficulty to show the relations for $e=1,2$.
Assume $e\geq 3$.
First we deal with the elements in the first row,
\begin{eqnarray*}
a_{1,j}^{(e+1)}
&=&
\sum_{s=1}^n a_{1,s}^{(e)} a_{s,j}=\sum_{s=1}^n U_{e-1}^{n-s}
U_e^{s-1} m^{s+j-n-1}
\binom{n-1}{s-1}\binom{s-1}{n-j}\\
&=&
\sum_{s=1}^n  U_{e-1}^{n-s} U_e^{s-1} m^{s+j-n-1}
 \binom{n-1}{j-1}\binom{j-1}{n-s}\\
&=&m^{j-1}
U_e^{n-1} \binom{n-1}{j-1}
\sum_{s=1}^n \left( \frac{U_{e-1}}{mU_e}\right)^{n-s} \binom{j-1}{n-s}\\
&=&
m^{j-1}
U_e^{n-1} \binom{n-1}{j-1} \left(1+\frac{U_{e-1}}{mU_e}\right)^{j-1}=
U_{e}^{n-j} U_{e+1}^{j-1} \binom{n-1}{j-1}.
\end{eqnarray*}
Now we prove the result for the elements in the
first column.
\begin{eqnarray*}
a_{i,1}^{(e+1)}
&=&
\sum_{s=1}^n a_{i,s} a_{s,1}^{(e)}=\sum_{s=1}^n
m^{i+s-n-1}
\binom{i-1}{n-s} U_{e-1}^{n-s} U_e^{s-1} \\
&=&
m^{i-1}
U_e^{n-1} \sum_{s=1}^n \left(\frac{U_{e-1}}{mU_e}
\right)^{n-s}\binom{i-1}{n-s}\\
&=&m^{i-1} U_e^{n-1}
\left(1+\frac{U_{e-1}}{mU_e}\right)^{i-1}
=
U_e^{n-i} U_{e+1}^{i-1}.
\end{eqnarray*}
Using \refeq{first_row}, we get
\begin{eqnarray*}
\sum_{j\geq 1} a_{1,j}^{(e)} y^{j-1}
&=&
\sum_{j\geq 1} \binom{n-1}{j-1} U_{e-1}^{n-j} U_e^{j-1}  y^{j-1}\\
&=&
\sum_{s\geq 0} \binom{n-1}{s} U_{e-1}^{(n-1)-s} (y U_e)^{s}\\
&=&
(U_{e-1}+yU_e)^{n-1}
\end{eqnarray*}
Using \refeq{first_eq_gen} and the fact that
$ (U_{e-1}-U_e x) \sum_{i\geq 1}
U_{e-1}^{n-i} U_e ^{i-1} x^{i-1}=U_{e-1}^n$ and $U_{e-1} a_{1,1}^{(e)}=U_{e-1}^n$,
 we deduce the result.
\end{proof}

%\newpage

%\section{An Inverse}

The inverse of $T_n(m)$ is not difficult to find. We have
\begin{thm}
The inverse of
$T_n(m)=\left(m^{i+j-n-1}\binom{i-1}{n-j}\right)_{i,j}$ is
\[
 T_n(m)^{-1}=
\left((-1)^{n+i+j+1} m^{n+1-i-j}\binom{n-i}{j-1}\right)_{i,j}.
\]
\end{thm}
\begin{proof}
The $(i,j)$ entry in $A^{-1} A$ is
\[
\begin{split}
& \sum_{s=1}^n (-1)^{n+i+s+1} m^{n+1-i-s}\binom{n-i}{s-1}
m^{s+j-n-1}\binom{s-1}{n-j} \\
& =m^{j-i} \sum_{s=1}^n (-1)^{n+i+s+1}
\frac{(n-i)!}{(s-1)!(n-i-s+1)!}
\frac{(s-1)!}{(n-j)!(s+j-n-1)!}\\
& =m^{j-i}  \sum_{s=1}^n (-1)^{n+i+s+1}
\frac{(n-i)!}{(n-j)!(j-i)!}
\frac{(j-i)!}{(n-i-s+1)!(s+j-n-1)!}%\\
\end{split}
\]
\[
\begin{split}
& =m^{j-i} \binom{n-i}{n-j}\sum_{s=1}^{n} (-1)^{n+i+s+1}
\binom{j-i}{n-i-s+1}\\
& =m^{j-i} \binom{n-i}{n-j}
\sum_{k=0}^{j-i} (-1)^{k}  \binom{j-i}{k},
\end{split}
\]
which is $0$, unless $i=j$, in which case it is $1$.
\end{proof}

\section{Powers of $T_n(m)$ modulo $p$}

Let $m\in {\bf Z}$.
Using the recurrence among the entries of $T_n(m)$,
and reasoning as in \cite{PS}, we prove the following
\begin{thm}
\label{the1_U}
If $e$ is the least integer (entry point) such that
$U_e \equiv 0 \pmod{p}$, then
\[
\begin{split}
T_{2k}(m)^e &\equiv (-1)^{(k+1)e} U_{e-1}I_{2k}\pmod p\\
T_{2k+1}(m)^e &\equiv (-1)^{ke}I_{2k+1} \pmod p.
\end{split}
\]
Moreover, $T_n(m)^{4e}\equiv I_n\pmod{p}$.
Furthermore, considering the parity of $e$, we have
\[
T_n(m)^{2e}\equiv I_n\pmod p\ \text{if $e$ even}
\]
and if $e$ odd
\[
\begin{split}
T_n(m)^{2e}\equiv &\, r^{n-1} I_n\pmod p\ \text{if $e\equiv 3\pmod 4$}  \\
T_n(m)^{2e} \equiv  &\,  (-r)^{n-1} I_n\pmod p\ \text{if $e\equiv 1\pmod 4$},
\end{split}
\]
where $r\equiv \frac{U_{(e+1)/2}}{U_{(e-1)/2}}\pmod p$, so
$r^2\equiv -1\pmod p$.
\end{thm}
\begin{proof}
Using \refeq{gen_fib}, if $U_e\equiv 0\pmod p$, then
\[
U_{e-1}a_{i,j}^{(e)}\equiv U_{e+1}a_{i-1,j-1}^{(e)}.
\]
Since $p$ divides neither $U_{e-1}$ nor $U_{e+1}$
(otherwise it would divide
$U_1=1$), we get
\[
\begin{split}
& a_{i,j}\equiv 0\pmod p,\ \text{if}\ i\not= j,\\
& a_{i,i}^{(e)}\equiv a_{i-1,i-1}^{(e)}
\equiv\cdots \equiv a_{1,1}^{(e)}\equiv U_{e-1}^{n-1}\pmod p.
\end{split}
\]
Therefore
\[
T_n(m)^{e}\equiv U_{e-1}^{n-1} I_n\pmod p.
\]
Using Cassini's identity $U_{l-1} U_{l+1}-U_l^2=(-1)^l$, for
$l=e$,
we get, if $n=2k$,
\[
%\begin{split}
U_{e-1}^{n-1}
=U_{e-1}^{2k-1}\equiv \left( U_{e-1}^2\right)^k U_{e-1}^{-1}
\equiv  (-1)^{ke} U_{e-1}^{-1}\equiv (-1)^{(k+1)e} U_{e-1}\pmod p,
%\end{split}
\]
since $U_{e-1}^2\equiv U_{e+1}^2\equiv (-1)^e \pmod p$.
If $n=2k+1$, then
\[
%\begin{split}
U_{e-1}^{n-1}
= U_{e-1}^{2k}\equiv \left( U_{e-1}^2\right)^k
\equiv  (-1)^{ke} \pmod p.
%\end{split}
\]
The previous two congruences replaced in
$T_n(m)^{e}\equiv U_{e-1}^{n-1} I_n\pmod p$, proves the first claim.

Lemma 3.4 of \cite{Li} implies
\[
\begin{split}
U_{e-1}&\equiv (-1)^{\frac{e-2}{2}}\ \text{if $e$ even}\\
U_{e-1}&\equiv r(-1)^{\frac{e-3}{2}}, r^2\equiv -1\pmod p,\
\text{if $e$ odd}.
\end{split}
\]
The residue $r$ in the previous relation is just
$r\equiv \frac{U_{(e+1)/2}}{U_{(e-1)/2}}\pmod p$.
Thus, if $e$ even, then $U_{e-1}^2\equiv 1\pmod p$, so
\[
T_n(m)^{2e}\equiv I_n\pmod p,
\]
for any $n$. The remaining cases are similar.
\end{proof}

Similarly, we can prove
\begin{thm}
\label{the2_U}
1.  If  $p\,|U_{p-1}$,
 then $T_n(m)^{p-1} \equiv I_n \pmod p$.\\
2. If  $p\,|U_{p+1}$,
then $T_{2k+1}(m)^{p+1} \equiv I_{2k+1} \pmod p$ and
$T_{2k}(m)^{p+1} \equiv -I_{2k} \pmod p$.
\end{thm}
A consequence of Theorem 1 of \cite{A} is
\begin{lem}
\label{lemma-Ando}
For a prime $p$ which divides $f(x)=x^2-mx-1$ for an integer $x$,
the sequence $\{U_e\}_e$ has a period $p-1\pmod p$, provided $p$ is
not a divisor of $D=m^2+4$.
\end{lem}
Our final result is
\begin{thm}
Let $p$ be a prime divisor of $x^2-mx-1$, for some integer $x$ and
$\gcd(p, m^2+4)=1$. Then,
 $T_n(m)^{p-1}\equiv I_n\pmod p$.
\end{thm}
\begin{proof}
Straightforward, using Lemma \refth{lemma-Ando} and Theorem
\refth{general_fibonacci} or Theorem \refth{the1_U}.
\end{proof}

%There are some very interesting applications to the previous theorem.
%We state a few in
%\begin{thm}
%The sequence
%\end{thm}

\section{Further Research}
We observed that netted matrices defined using
second/third-order recurrences (we call these $2$ or
$3-netted$ matrices) preserve a third-order
recurrence among the entries of their powers. The natural question arising
is: {\em what is the degree of the recurrence (if it exists - we conjecture
that it does) for higher powers  of a $4$-netted, $5$-netted, etc.,
matrix?}.

The spectral properties of $T_n(m)$ is another topic of future research.
Let $U_e=U_e(m)$ be the general Pell or Fibonacci sequence. It is known that
$U_e=\frac{\alpha^e-\beta^e}{\alpha-\beta}$, where
$\alpha=\frac{m+\sqrt{m^2+4}}{2},\, \beta=\frac{m-\sqrt{m^2+4}}{2}$.
We associate
 the general Lucas sequence $V_e$ satisfying the same recurrence as $U_e$, with initial
 conditions $V_0=2$, $V_1=m$. Thus $V_e=\alpha^e+\beta^e$.
 We conjecture
 \begin{conj}
The characteristic polynomial $p_n(x)$ of $T_n(m)$ is
\[
\begin{split}
&p_1(x)=1-x;\
p_{4k+1}=(1+ V_{4k-2}x+x^2)(1-V_{4k}x+x^2)p_{4k-3}(x)\\
&p_3(x)=-(1+x)(1-V_2 x+x^2);\
p_{4k+3}=(1+V_{4k}x+x^2)(1-V_{4k+2}x+x^2)p_{4k-1}(x)\\
&p_2(x)=-1- V_1 x + x^2;\
p_{4k+2}=(-1+V_{4k-1}x+x^2)(-1- V_{4k+1}x+x^2)p_{4k-2}(x)\\
&p_0(x)\stackrel{def}{=}1;\
p_{4k+4}=(-1+V_{4k+1}x+x^2)(-1- V_{4k+3}x+x^2)p_{4k}(x).
\end{split}
\]
 \end{conj}
We checked the conjecture up to dimension $100\times 100$.


\begin{thebibliography}{99}

\bibitem{A} S. Ando, {\em On the Period of Sequences Modulo a
Prime Satisfying a Second Order Recurrence}, Applications of
Fibonacci Numbers, Vol. 7, 1998, pp. 17-22.

\bibitem{D} L.E. Dickson,
History of the Theory of Numbers, Vol. 1, Ch. XVII,
Chelsea Publishing Co., 1971.

\bibitem{GKP} R.L. Graham, D.E. Knuth, O. Patashnik,
Concrete Mathematics,
Adison-Wesley Publishing Company, 1989.

\bibitem{Li} H.-C. Li, {\em On Second-Order Linear Recurrence
Sequences: Wall and Wyler Revisited}, Fibonacci Quarterly,
Nov. 1999, pp. 342-349.

\bibitem{PS} R. Peele, P. St\u anic\u a,
{\em Matrix Powers of Column-Justified Pascal Tri\-angles and
Fi\-bo\-na\-cci
Se\-quen\-ces}, to appear %in Fi\-bo\-na\-cci Quar\-ter\-ly
(available at {\em http://sciences.aum.edu/\~\,stanpan}).

\bibitem{AeqB} M. Petkovsek, H. Wilf, D. Zeilberger, $A=B$,
A.K. Peters, Ltd., Wellesley, Massachussets, 1997.

\end{thebibliography}
\end{document}